\documentclass[11pt]{article} 
\newcommand{\version}{May 25, 2010}

\usepackage{amsthm,amsfonts,amsmath,amssymb}
\usepackage{graphicx,color}
\swapnumbers
\theoremstyle{plain}

\newtheorem{thm}{THEOREM}[section]
\newtheorem{lm}[thm]{LEMMA}

\theoremstyle{definition}
\newtheorem{defi}[thm]{DEFINITION}
\theoremstyle{remark}
\newcommand{\eps}{\ensuremath{\varepsilon}}
\newcommand{\upchi}{\raise1pt\hbox{$\chi$}}
\newcommand{\R}{{\mathord{\mathbb R}}}

\newcommand{\N}{{\mathord{\mathbb N}}}

\renewcommand{\|}{{\Vert}}
\numberwithin{equation}{section}  
\def\dd{{\rm d}}

\begin{document}

\markboth{\scriptsize{CU \version}}{\scriptsize{CU \version}}

\title{\bf{ Localization, Smoothness, and Convergence to Equilibrium for
 a Thin Film Equation}}
\author{\vspace{5pt} Eric A. Carlen$^{1}$ and S\"{u}leyman Ulusoy$^{2}$
 \\
\vspace{5pt}\small{1 Department of Mathematics, Rutgers University, NJ, USA,}\\
\small{2 CSCAMM, University of Maryland, College Park, MD, USA }
\\}
\date{\version}
\maketitle \footnotetext [1]{Work partially supported by U.S.
National Science Foundation
grant DMS DMS  0901632.    \\
\maketitle \footnotetext [2]Work partially supported by U.S.
National Science Foundation\\
grants DMS 0707949, DMS1008397 and FRG0757227.  \\
\copyright\, 2010 by the authors. This paper may be reproduced, in
its entirety, for non-commercial purposes.}

\begin{abstract}
We investigate the long-time behavior of weak solutions to the thin-film type equation 
$$v_t =(xv - vv_{xxx})_x\ ,$$  
which arises in the Hele-Shaw problem. We estimate the rate of convergence of solutions  to
the Smyth-Hill equilibrium solution, which has the form $\frac{1}{24}(C^2-x^2)^2_+$,  in the norm
$$|\!|\!| f |\!|\!|_{m,1}^2 = \int_{\R}(1+ |x|^{2m})|f(x)|^2\dd x +  \int_{\R}|f_x(x)|^2\dd x\ .$$
We obtain exponential convergence in the $|\!|\!| \cdot |\!|\!|_{m,1}$ norm for all $m$ with $1\leq m< 2$, thus obtaining rates of convergence in norms measuring both smoothness and localization. The localization is the main novelty, and in fact, we show that there is a close
connection between the localization bounds and the smoothness bounds: Convergence of second moments implies convergence in the
$H^1$ Sobolev norm.   We then use methods of optimal mass transportation to obtain the convergence of the required moments.
We also use such methods to construct an appropriate class of weak solutions for which all of the estimates on which our convergence analysis depends may be rigorously derived. Though  our main results on convergence can be stated without reference to optimal mass transportation,
essential use of this theory is made throughout our analysis.

\end{abstract}

\bigskip
\centerline{Key words: thin-film equation, Wasserstein distance, gradient flow, Euler-Lagrange equation }
\centerline{Mathematics Subject Classification Numbers: 35A15, 35B40,
35K25, 35K45, 35K55, 35K65 }

\newpage

\section{Introduction} \label{intro}

\subsection{The primary Lyapunov functional}

The following
one-dimensional fourth-order nonlinear degenerate parabolic
equation

\begin{equation}\label{biri}
u_t = -\left(u u_{xxx} \right)_x, \qquad  x  \in  \mathbb{R}, t>0,\
\end{equation}
with
\begin{equation}\label{ini}
u(x,0) = u_0(x) \geq 0, \qquad x \in \mathbb{R},
\end{equation} arises as the particular case of
the  thin-film equation in the the Hele-Shaw setting \cite{B,M,O}.

Equation (\ref{biri}) has a well-known scale invariance  and corresponding self-similar solutions \cite{CT,SH}: Given a 
solution $u(x,t)$ of (\ref{biri}), define $v(x,t)$ by
\begin{equation}\label{v}
v(x,t) := a(t) u(a(t)x, b(t)),
\end{equation}
where
$a(t) = e^t $ and $ b(t) = (e^{5t} -1)/5$.
Then  $v(x,t)$ satisfies:
\begin{equation}\label{veq}
v_t = (xv - v v_{xxx})_x, \qquad x \in \mathbb{R}, t>0,
\end{equation}
\begin{equation}\label{vic}
v(x,0) = u_0(x), \qquad x \in \mathbb{R}.
\end{equation}

Both equations (\ref{biri}) and (\ref{veq}) are conservation laws, so that the total mass is conserved. That is, in the case of (\ref{veq}),
$$M = \int_{\R}  v(x,t){\rm d}x $$
is conserved. 

Also, both equations describe a gradient flow \cite{CT}: For (\ref{veq}), define the energy functional $E[v]$ where
$$ E[v] = \frac{1}{2} \int_{\mathbb{R}} \left(  v_x^2(x) + x^2 v(x) \right) \, \dd x.$$
Then, (\ref{veq}) can be written as
\begin{equation}\label{vgrf}
 v_t = \left (  v \left ( \frac{\delta E[v]}{\delta v} \right )_x   \right)_x\ .
 \end{equation}
 
It follows from (\ref{vgrf}) that for sufficiently regular solutions of (\ref{veq}),
$$\frac{{\rm d}}{{\rm d}t}E[v] = -\int_\R v \left[  \left ( \frac{\delta E[v]}{\delta v} \right )_x   \right]^2\dd x\ .$$
Thus, as long as the solution $v$ is non-negative,  $E[v]$ will be monotone decreasing.  

The physical interpretation of $u$ as the height of a thin film  \cite{B,M,O} suggests that $u$, and hence $v$, must be non-negative for all physically meaningful solutions. However, there is no maximum principle  to provide {\it a-priori} assurance of non-negativity for this fourth order equation, and it is an open question whether  certain solutions, necessarily non-physical, can become negative.  
In what follows, we shall construct and analyze a particular class of non-negative solutions.

In the rest of this introduction, we restrict our attention to non-negative solutions for which  $E[v]$ is monotone decreasing.

It is easy to determine the minimizers of $E[v]$, each of which is of course a stea\dd y state solution
of (\ref{veq}), and thus determines a self-similar solution of (\ref{biri}): These minimizing steady states \cite{SH}
are the one parameter family of functions
\begin{equation}\label{ssv}
v^{(\infty)}(x) := \frac{1}{24}\left( C^2 - x^2   \right)_{+}^2,
\end{equation}
where the subscripted $+$ indicates the positive part,  and the parameter  $C$ determines the total mass of
$v^{(\infty)}$. Because mass is conserved,  we suppress $C$ in our notation.

Define the {\em relative energy} $E[v|v^{(\infty)}]$ by
\begin{equation}\label{relen}
E[v|v^{(\infty)}] = E[v] - E[v^{(\infty)}]\ .
\end{equation}
It turns out that this is a good Lyapunov functional for the steady states $v^{(\infty)}$:

\begin{lm}\label{energyalph} Let $v$ have the same mass as $v^{(\infty)}$. Then:
\begin{equation}\label{h1}
E[v|v^{(\infty)}] \geq  \frac{1}{2} \int_{\mathbb{R}} |v_x - v_x^{(\infty)}|^2 \, \dd x + 
\frac{1}{3} \int_{\{v^{(\infty)}=0\}}|x|^2v \, \dd x\ .
\end{equation}
\end{lm}

\noindent{\bf Proof:}  One finds, after one integration by parts, and using the fact that $v$ and $v^{(\infty)}$ have the same mass.
\begin{eqnarray}
E[v|v^{(\infty)}] &-&  \frac{1}{2} \int_{\mathbb{R}} |v_x - v_x^{(\infty)}|^2 \, \dd x\nonumber\\
&= &
-\frac{1}{2} \int_{\mathbb{R}}|x|^2v^{(\infty)} \, \dd x + \frac{1}{6}C^2 \int_{\mathbb{R}}v^{(\infty)} \, \dd x
-\int_{\mathbb{R}}(v_x^{(\infty)})^2 \, \dd x\nonumber\\
&-& \int_{\{v^{(\infty)}=0\}}\left(\frac{1}{6}C^2-  \frac{1}{2}|x|^2\right)v \, \dd x\ .\nonumber
\end{eqnarray}
By direct computation, the three integrals on the second line cancel exactly, and on ${\{v^{(\infty)}=0\}}$, $|x|^2/2 - C^2/6\geq |x|^2/3$. \qed
\medskip

Thus,  $E[v|v^{(\infty)}]$ is monotone decreasing along the evolution described by (\ref{veq}), and by (\ref{h1}), whenever 
$E[v|v^{(\infty)}]$ is small, $v$ is close to $v^{(\infty)}$, which is what we mean by saying that  $E[v|v^{(\infty)}]$ is a good Lyapunov functional.
Let us express this quantitatively. We shall use the following norms:

\begin{defi} For any smooth compactly supported function $f$ on $R$, and any $m\geq 0$, define 
\begin{equation}\label{normdef}
|\!|\!| f |\!|\!|_{m,1}^2 = \int_{\R}(1+ |x|^{2m})|f(x)|^2\dd x +  \int_{\R}|f_x(x)|^2\dd x\ .
\end{equation}
We extend then the norm $|\!|\!| \cdot |\!|\!|_{m,1}$ to its natural Hilbert space domain. 
\end{defi}

\begin{lm}\label{energyalph2} Let $v\geq 0$ have the same mass $M$ as $v^{(\infty)}$.  Then:
\begin{equation}\label{infbnd}
\|v\|_\infty \leq  \frac{5}{3}M^{3/5}(2E[v])^{1/5}\ ,
\end{equation}
and there exist explicitly computable constants $K_1$ and $K_2$ depending only on $M$ such that 
\begin{equation}\label{h1a}
|\!|\!| v- v^{(\infty)} |\!|\!|_{1,1}^2 \leq   K_1 (E[v|v^{(\infty)}])^2 +  K_2 E[v|v^{(\infty)}]\ .
 \end{equation}
\end{lm}

\noindent{\bf Proof:}  For any $x > y$,
${\displaystyle v(x) = v(y) + \int_x^y v_x(z)\dd z \leq v(y) + \sqrt{x-y}\|v_x\|_2}$.
Now for any $a>0$,  integrating in $y$ over $[x-a,x]$,
${\displaystyle av(x) \leq \int_{x-a}^x v(y)\dd y + \frac{2}{3}a^{3/2}\|v_x\|_2}$.  Optimizing in $a$, one obtans
(\ref{infbnd}). Next, 
\begin{eqnarray}\label{inin1}
\int_{\{|x|\geq C\}}(1+x^2)(v- v^{(\infty)})^2 \dd x &\leq& \|v\|_\infty  \int_{\{|x|\geq C\}}(1+x^2)v\dd x\nonumber\\
&\leq& \frac{C^2+1}{C^2}\|v\|_\infty   \int_{\{|x|\geq C\}}x^2v\dd x\nonumber\\
&\leq& \frac{C^2+1}{C^2}\|v\|_\infty 3E[v|v^{(\infty)}]\ ,
\end{eqnarray}
where we have used  (\ref{h1}) in the last line.
Next, since $v$ and $v^{(\infty)}$ have the same mass,
$$\left|\frac{1}{2C} \int_{\{|x| \leq C\}} (v- v^{(\infty)})\dd  x\right| =  \frac{1}{2C} \int_{\{|x|\geq C\}}v\dd x \leq \frac{1}{2C^3} \int_{\{|x|\geq C\}}|x|^2v\dd x \ .$$
Hence there exists an $x_0\in [-C,C]$ such that  
${\displaystyle |v(x_0) - v^{(\infty)}(x_0)| \leq  \frac{1}{2C^3} \int_{\{|x|\geq C\}}|x|^2v\dd x}$.
Thus, for any $x\in [-C,C]$,  
\begin{eqnarray}
|v(x) - v^{(\infty)}(x)| &\leq&  |v(x_0) - v^{(\infty)}(x_0)| + \sqrt{|x-x_0|}\|v_x - v_x^{(\infty)}\|_2\nonumber\\
&\leq& \frac{3}{2C^3}E[v|v^{(\infty)}]  + 2\sqrt{CE[v|v^{(\infty)}] }\nonumber
\end{eqnarray}
where in the last line we have used (\ref{h1}) twice. It follows that
\begin{equation}\label{inin2}
\int_{\{|x|\leq C\}}(1+x^2)(v- v^{(\infty)})^2 \dd x \leq  (1+C^2)\left[ \frac{9}{C^2}(E[v|v^{(\infty)}])^2 + 8C^2E[v|v^{(\infty)}]\right]\ .
\end{equation}
Combining  (\ref{infbnd}),
(\ref{inin1}) and (\ref{inin2} ), and recalling that $M$ determines $C$, we obtain the result, and see that explicit values of $K_1$ and $K_2$
may be written down. 
\qed

Hence any result showing that $E[v|v^{(\infty)}]$ converges to zero along solutions of (\ref{veq}) shows that $v$ converges to $v^{(\infty)}$ in the 
$|\!|\!| \cdot |\!|\!|_{1,1}$ norm.  Such results have been proved in \cite{CU} and more  recently  in \cite{MMS}, with an optimal exponential rate. 
However, in these papers, a weaker form of (\ref{h1}) was used, without the second term on the right. Consequently, these papers only deduced
the convergence of $v$ to $v^{(\infty)}$ in the $|\!|\!| \cdot |\!|\!|_{0,1}$ norm.

At first sight, it may seem a trivial matter to go from convergence in the $|\!|\!| \cdot |\!|\!|_{0,1}$ norm to convergence in the
$|\!|\!| \cdot |\!|\!|_{m,1}$ norm for higher values of $m$:
One might guess that since the steady state $v^{(\infty)}$ is supported in $[-C,C]$, it should be easy to control the evolution of higher moments
$ M_{2m}(v) := \int_{\mathbb{R}}|x|^{2m}v(x,t) \, \dd x$ of solutions $v$ of (\ref{veq}). Upon careful consideration, this turns out not to be the case. We shall present a non-trivial argument to show that  $ M_4(v(\cdot,t))$ is bounded uniformly in $t$ in terms of its initial value, but we do not know if this is even true for
 moments higher than the fourth.  So while one might expect to be able to prove strong localization estimates for solutions of (\ref{veq}), localization and moment bounds turn out to be somewhat subtle. Using our uniform bound on  $M_4(v)$, we shall deduce an exponential rate of convergence of 
 $v$ to  $v^{(\infty)}$ in the $|\!|\!| \cdot |\!|\!|_{m,1}$ norm for all $1 \leq m < 2$.
 
Moreover, we shall give a new proof of the fact that  $E[v|v^{(\infty)}]$ converges to zero at an exponential rate. In this new proof, moment estimates play the crucial role:  The rate at which  $E[v|v^{(\infty)}]$ converges to zero is controlled by the rate at which 
$\int_{\mathbb{R}}|x|^{2}v\dd x$ converges to $\int_{\mathbb{R}}|x|^{2}v^{(\infty)}(x) \dd x$.  This is somewhat remarkable: We are using a simple functional involving no derivatives of $v$ to control one that does involve derivatives of $v$. There turns out to be close interplay be between smoothness and localization for solutions of (\ref{veq}), and one point of this paper is to explain this interplay, and show how it may be used.

 \subsection{Convergence of second moments and convergence of the energy}

Let $v(x,t)$ be a solution of (\ref{veq}), and let  $v^{(\infty)}$ be the stationary solution of the same total mass.  Then the {\em relative second moment}
\begin{equation}\label{relpen}
\alpha[v|v^{(\infty)}] =  \alpha[v(\cdot,t)] - \alpha[v^{(\infty)}] \qquad{\rm where}\qquad
 \alpha[v] = \frac{1}{2} \int_{\mathbb{R}} x^2 v(x) \, \dd x\ ,
\end{equation}
and the {\em relative surface energy}
\begin{equation}\label{relien}
\beta[v|v^{(\infty)}] =  \beta[v(\cdot,t)] - \beta[v^{(\infty)}] \qquad{\rm where}\qquad
 \beta[v] = \frac{1}{2} \int_{\mathbb{R}} v_x^2 \, \dd x\, \dd x\ .
\end{equation}

Evidently, the three quantities $E[v|v^{(\infty)}]$,  $\alpha[v|v^{(\infty)}]$ and $\beta[v|v^{(\infty)}]$ are related by
\begin{equation}\label{rela}
E[v|v^{(\infty)}] =  \alpha[v|v^{(\infty)}] +  \beta[v|v^{(\infty)}]\ .
\end{equation}
It follows that if one can show that any two of these converge to zero, so does the third. As we have seen above, 
 if one can show that $\lim_{t\to\infty} E[v|v^{(\infty)}] =0$, one concludes as well that \hfill\break $\lim_{t\to\infty} |\!|\!| v - v^{(\infty)} |\!|\!|_{1,1}= 0$,
 from which it certainly follws that 
$\lim_{t\to\infty} \alpha[v|v^{(\infty)}] =0$.

What is perhaps more surprising is that if one can show that $\lim_{t\to\infty} \alpha[v|v^{(\infty)}] =0$, one can also deduce as a direct consequence
that $\lim_{t\to\infty} E[v|v^{(\infty)}] =0$, and moreover, one can estimate the rate of convergence in the latter limit in terms of the former limit. 
Let us explain how this works, first at the level of formal calculation.

We easily compute that \cite{CU}
\begin{equation}
 \frac{d}{dt} \alpha[v(\cdot,t)|v^{(\infty)}] = -2 \alpha[v(\cdot,t)|v^{(\infty)}]  + 3 \beta[v(\cdot,t)|v^{(\infty)}]\ .
 \end{equation}
From this and (\ref{rela}) we get
\begin{equation}\label{ale}
\frac{d}{dt} \alpha[v(\cdot,t)|v^{(\infty)}] = -5 \alpha[v(\cdot,t)|v^{(\infty)}] + 3 E[v(\cdot,t)|v^{(\infty)}].
\end{equation}
Now, for $T>1$ let us integrate both sides of (\ref{ale}) from $T-1$ to $T$ to obtain
\begin{equation}\label{expd}
\alpha(v(\cdot,T)|v^{(\infty)}) - \alpha(v(\cdot,T-1)|v^{(\infty)}) + 5 \int_{T-1}^{T} \alpha[v(\cdot,t)|v^{(\infty)}] \, dt \geq 3  E[v(\cdot,T)|v^{(\infty)}],
\end{equation}
since $E[v(\cdot,t)|v^{(\infty)}]$ is monotone decreasing. 
Now suppose we have an estimate of the form
\begin{equation}\label{expd2}
|\alpha(v(\cdot,t)|v^{(\infty)})| \leq K e^{-\lambda t}\ .
\end{equation}
Using this in (\ref{expd}) yields
\begin{equation}\label{expd3}
  E[v(\cdot,T)|v^{(\infty)}] \leq \frac{K}{3}\left(1 + e^{\lambda} + \frac{5}{\lambda}e^{\lambda}\right)e^{-\lambda T}\ .
\end{equation}
In the next subsection, we explain how we shall obtain a rate of convergence estimate for  $\alpha[v(\cdot,t)|v^{(\infty)}]$.

\subsection{The second Lyapunov functional}

Carrillo and Toscani have made the remarkable discovery \cite{CT} that the equation (\ref{veq}) possesses a second Lyapunov functional: Define
the {\em entropy} $H[v]$ by 
$$ H[v] = \int_{\mathbb{R}} \left( \frac{x^2}{2} v(x) + 2 \sqrt{\frac{2}{3}} v^{3/2}(x) \right) \, \dd x.$$
and then  the relative entropy by
$H[v|v^{(\infty)}] = H[v] - H[v^{(\infty)}]$.
We remark that by (\ref{infbnd}),
\begin{equation}\label{impl}
E[v] < \infty \Rightarrow H[v] < \infty\ .
\end{equation}
It is easy to see that  $v^{(\infty)}$ minimizes $H[v]$ among all non-negative integrable functions $v$ with the same mass as 
$v^{(\infty)}$, and hence $ H[v|v^{(\infty)}] $ is non-negative.  In fact, as shown by Otto \cite{otto},
\begin{equation}\label{pkl}
H[v|v^{(\infty)}] \geq \left(\int_\R |v - v^{(\infty)}|\dd x\right)^2\ .
\end{equation}

The entropy functional $H$ arises in the theory of the porous medium equation. There is a particular slow-diffusion case of the porous medium equation for 
which the Smyth-Hill densities $v^{(\infty)}$ are also steady state solutions. This equation, which can be written in the gradient flow form
\begin{equation}\label{pmgrf}
 v_t = \left (  v \left ( \frac{\delta H[v]}{\delta v} \right )_x   \right)_x\ ,
\end{equation}
is second-order parabolic. For it, the maximum principle applies and provides both positivity and uniqueness. Hence $H[v|v^{(\infty)}]$ is a Lyapunov functional for the equation (\ref{pmgrf}), and on account of (\ref{pkl}), it is a good one. 

The remarkable discovery of Carrillo and Toscani is that $H[v|v^{(\infty)}]$ is also a Lyapunov functional for the thin film equation (\ref{veq}), and indeed they even show that for strong solutions of (\ref{veq}) and $T>S$,
\begin{equation}\label{strongL}
H[v(\cdot,T)|v^{(\infty)}] + 2\int_S^T H [v(\cdot, t)|v^{(\infty)}] \dd t \leq H[v(\cdot,S)|v^{(\infty)}] \ .
\end{equation}
This has the immediate consequence that
\begin{equation}\label{strongL2}
H[v(\cdot,t)|v^{(\infty)}]  \leq e^{-2t}H[v(\cdot,0)|v^{(\infty)}] \ .
\end{equation}

The fact that (\ref{strongL}) holds for solutions $v$ of the slow diffusion equation (\ref{pmgrf}) was discovered by Otto \cite{otto} using methods from the theory of optimal mass transportation, and in particular, the notion of {\em displacement convexity}, as we shall explain in more detail in the next section.

The fact that (\ref{strongL}) also holds for solutions of (\ref{veq}) is far from obvious, but was proved by Carrillo and Toscani 
using the fact that 
the  equation (\ref{veq}) can be written as 
$$ v_t = -\left( \Phi(v) \left[ \frac{x^2}{2}+h(v)\right]_{xx}\right)_{xx} + \left( v \left[ \frac{x^2}{2} + h(v) \right]_x \right)_x\ ,$$
with $h(v) = \sqrt{6}v^{1/2}$ and $\Phi(v) = v h'(v)$.  Combining  (\ref{strongL2}) and (\ref{pkl}), they then deduced
\begin{equation}\label{strongL3}
\int_\R |v - v^{(\infty)}|\dd x \leq e^{-t}\sqrt{H[v(\cdot,0)|v^{(\infty)}]}\ ,
\end{equation}
and raised the question of proving that the energy $E[v(\cdot,t)|v^{(\infty)}]$ converges to zero.

As explained in the previous subsection, to do this, it suffices to prove that $\alpha[v(\cdot,t)|v^{(\infty)}]$ converges to zero.
We shall do this by proving in the second section of this paper an inequality relating for the entropy
$H[v|v^{(\infty)}] $ and $\alpha[v(\cdot,t)|v^{(\infty)}]$.  We shall prove:

\begin{lm}\label{entsecmo} For any non-negative integrable  function $v$ on $\R$ with a finite second moment, let $v^{(\infty)}$ be the Smyth-Hill density with the same mass. 
Then 
$$\alpha[v|v^{(\infty)}] \leq   2\sqrt{\alpha[v^{(\infty)}]}  \sqrt{  H[v|v^{(\infty)}]  } +  H[v|v^{(\infty)}] \ .$$
\end{lm}

Granted this lemma, the bound (\ref{strongL2}) now gives us
$$\alpha[v(\cdot,t)|v^{(\infty)}] \leq   e^{-t}2\sqrt{\alpha[v^{(\infty)}]}  \sqrt{  H[v(\cdot,0)|v^{(\infty)}]  } + e^{-2t}H[v(\cdot,0)|v^{(\infty)}] \ ,$$
and then from (\ref{expd3}) we have
$E[v(\cdot,t)|v^{(\infty)}] \leq Ce^{-t}$ for an explicit constant $C$. 
Finally, by Lemma~\ref{energyalph2}, we obtain exponential convergence in the $|\!|\!| \cdot  |\!|\!|_{1,1}$ norm.
This outlines the general scheme of our strategy. 

We remark that in \cite{MMS} it is shown that for a certain class of solutions, $E[v(\cdot,t)|v^{(\infty)}]$ decays like $e^{-2t}$, which is twice the rate implied by
our result that $\alpha[v(\cdot,t)|v^{(\infty)}]$ decays like $e^{-t}$. However, our construction is somewhat less delicate, and the information we obtain on moments
and localization is new.

\subsection{Properly dissipative weak solutions}

The theory of the thin film equation is not yet in a well developed state. Basic issues of existence of strong solutions and uniqueness remain open. For which class of solutions can the formal calculations above be made precise and rigorous? We now introduce such a class of weak solutions, called \emph{properly dissipative weak solutions}.  

The key estimates used in our convergence analysis depend on the  fact that 
$\alpha[v(\cdot,t)|v^{(\infty)}]$ satisfies (\ref{ale}) and that $H[v(\cdot,t)|v^{(\infty)}]$ satisfies (\ref{strongL}). Therefore, what we need is existence of weak solutions with these properties.

\begin{defi}\label{wsol}
A non-negative measurable function $v \in C^0(\mathbb{R} \times [0,\infty))$ such that for some fixed $M$, the {\em mass}, and each $t\geq 0$,
$$\int_\R v(x,t)\dd x = M \qquad{\rm and}\qquad E[v(\cdot, t)] < \infty\ ,$$
 is called a \emph{ weak solution} if
 for all $\zeta \in C_0^{\infty}(\mathbb{R} \times (0, \infty))$
\begin{equation}\label{int-eqn}
\iint_{\mathbb{R} \times (0, \infty)} \left(  -v\zeta_t - 3 ( v_x)^2 \zeta_{xx} - 2 v v_x \zeta_{xxx} \right) \, \dd t \, \dd x = 0.
\end{equation}
It is called a 
{\em properly dissipative weak solution} if moreover:

\begin{itemize}

\item For all $t>0$
\begin{equation}\label{H-dissp}
2\int_0^t H(v|v^{(\infty)}) \, dt + H[v(\cdot,t)|v^{(\infty)}] \leq H[v(\cdot,0)|v^{(\infty)}],
\end{equation}


\item For all $t\geq 1$,
\begin{equation}\label{alph-disp}
\alpha[v(\cdot,t)|v^{(\infty)}] - \alpha[v(\cdot,t-1)|v^{(\infty)}] + 5\int_{t-1}^t \alpha[v(\cdot,s)|v^{(\infty)}] \dd s \geq 3E[v(t)|v^{(\infty)}]\ .
\end{equation}

\end{itemize}

\end{defi}


Our main existence theorem for properly dissipative weak solutions is the following:

\begin{thm}\label{weakex}  Let $v_0$ be any non-negative integrable function on $\R$ such that $E[v_0] < \infty$ and such that
$$M_4(v_0) = \int_{\R}x^4v_0(x)\dd x < \infty\ .$$
Let $M$ be the total mass of $v_0$, and let $v^{(\infty)}$ denote the Smyth-Hill steady state with the same mass $M$. Then there exists a properly dissipative weak solution $v$ such that $v(x,0) = v_0(x)$.
Moreover, there is an explicit constant $K_3$ depending only on $M$ and $E[v_0]$
\begin{equation}\label{4prop}
\int_\R |x|^4 v(x,t)\dd x \leq  \int_\R |x|^4 v_0(x)\dd x + 
K_3 H[v_0|v^{(\infty)}]\ .
\end{equation}
\end{thm}

As noted in  in (\ref{impl}),  $H[v_0|v^{(\infty)}]$ is finite whenever $E[v_0]$ is finite. Theorem~\ref{weakex} is proved in the third section of the paper. We now state another of our main results:

\begin{thm}\label{weakpro}  Let $v_0$ be any non-negative integrable function on $\R$ such that $E(v_0) < \infty$ and such that
$M_4(v_0)<\infty$. 
Let $M$ be the total mass of $v_0$, and let $v^{(\infty)}$ denote the Smyth-Hill steady state with the same mass $M$. 
Then for any properly dissipative weak solution of (\ref{veq}) with $v(x,0) = v_0(x)$,
\begin{equation}\label{conv1}
|\!|\!| v - v^{(\infty)} |\!|\!|_{1,1}^2 \leq K_4 e^{-t}
\end{equation}
where $K_4$ is a positve constant depending only on $M$ and $E(v_0)$. Moreover, for all $1< p < 2$,
\begin{equation}\label{conv2}
|\!|\!| v - v^{(\infty)} |\!|\!|_{p,1}^2 \leq K_5 e^{-(2-p) t}
\end{equation}
where $K_5$ is a positive constants depending only on $p$,  $M$, $M_4(v_0)$ and $E(v_0)$. 
\end{thm}

This theorem is also proved in the third section.

\section{Mass transportation and the thin film equation}

We rely in an essential way on methods of optimal mass transportation to both construct and analyze our weak solutions. In this section we
 briefly recall the points that are essential here. See \cite{V} or \cite{AGS} for more information. 
 
For $M>0$, let ${\mathcal M}_M$ denote the set of non-negative Borel measure $\mu$ on $\R$ with $\mu(\R) = M$ and such that
$$\int_{\R}|x|^2\dd \mu(x)  < \infty\ .$$  For $\mu,\nu\in {\mathcal M}_M$, define $\Gamma(\mu,\nu)$ to be the set of all non-negative
Borel measures $\gamma$ on $\R^2$ such that for all Borel sets $A \subset \R$,
$$\gamma(A\times \R) = \mu(A)\qquad{\rm and}\qquad  \gamma(\R \times A) = \nu(A)\ .$$
The set $\Gamma(\mu,\nu)$ is the set of all {\em couplings} of $\mu$ and $\nu$. 

The $2$-Wasserstein distance between $\mu,\nu\in {\mathcal M}_M$,  ${\rm W}_2(\mu,\nu) $, is defined by
\begin{equation}\label{wrel}
 {\rm W}_2^2(\mu,\nu) = \inf_{\gamma \in \Gamma(\mu,\nu)} \left\{  \int_{\R^d \times \R^d}   |x-y|^2 \, \dd\gamma(x,y)    \right \}\ .
\end{equation}

By the Brenier-McCann Theorem, when $\mu$ and $\nu$ are absolutely continuous, the infimum is attained at a unique optimal coupling $\gamma_*$,
which is concentrated on the graph of the derivative $\varphi_x$ of a convex function $\varphi$:  For all Borel sets $A,B\subset \R$, 
\begin{equation}\label{coup}
\gamma(A\times B) = \mu(A\cap \varphi_x^{-1}(B))\ .
\end{equation}
 It follows that
\begin{equation}\label{was-def}
{\rm W}_2^2(\mu,\nu) = \int_{\R} |x- \varphi_x(x)|^2 \, \dd\mu(x).
\end{equation}
We write ${\rm W}_2^2(v,w)$ to denote
${\rm W}_2^2(v\dd x ,w\dd x)$, as is standard.

The fact that the optimal coupling is induced by a map $x\mapsto \varphi_x(x)$, which is $\mu$-almost everywhere invertible, yields an interpolation
between $\mu$ and $\nu$. Let $\varphi$ be the a convex function such that (\ref{coup}) defines an optimal coupling of $\mu$ and $\nu$. For $t\in [0,1]$, define the convex function $\varphi^{(t)}$ by
$$\varphi^{(t)}(x) = (1-t)\frac{x^2}{2} + t \varphi(x)\ .$$
Then $\varphi^{(t)}_x$ interpolates between the identity and $\varphi_x$, and we define and $\mu^{(t)} \in {\mathcal M}_M$ 
$$\mu^{(t)}(A) = \mu( (\varphi^{(t)}_x)^{-1}(A))\ .$$
The map $t\mapsto \mu^{(t)}$ is McCann's {\em  displacement interpolation} between $\mu$ and $\nu$, 

There is another way of expressing this that will be useful to us. Let $T:\R\to \R$ be measurable. Then $T\# \mu$, the {\em push-forward}
of $\mu$ under $T$, is the measure given by $(T\#\mu)(A) = \mu(T^{-1}(A))$. Thus, $\mu^{(t)} = (\varphi^{(t)})_x\# \mu$.

A functional  $\mathcal{G}$ on  ${\mathcal M}_M$ is said to be {\em $\lambda$ displacement convex}, \cite{Mc,V} if for all $\mu,\nu\in {\mathcal M}_M$,
the displacement interpolation of $\mu$ and $\nu$ satisfies
\begin{equation}
\lambda t(1-t){\rm W}_2^2(\mu,\nu) + \mathcal{G}(\mu^{(t)}) \leq (1-t) \mathcal{G}(\mu) + t\mathcal{G}(\nu),
\end{equation}
for all $0\leq t \leq 1$.  In this case,
\begin{equation}\label{abovetan}
\lambda {\rm W}_2^2(\mu,\nu) + \lim_{t\downarrow 0} \frac{ {\mathcal G}[\mu^{(t)}] - {\mathcal G}[\mu]}{t} \leq 
{\mathcal G}[\nu]- {\mathcal G}[\mu]\ .
\end{equation}
Then, if $\mu$ minimizes ${\mathcal G}$ so that the subgradient of ${\mathcal G}$ vanishes at $\mu$ (see\cite{AGS}),  (\ref{abovetan})
reduces to
\begin{equation}\label{abovetan2}
\lambda {\rm W}_2^2(\mu,\nu)  \leq 
{\mathcal G}[\nu]- {\mathcal G}[\mu]\ ,
\end{equation}
which is known as a 
{\em Talagrand inequality} for the  functional ${\mathcal G}$.

This is of immediate relevance to the entropy functional $H[v]$, regarded as a functional on ${\mathcal M}_M$ in the obvious way,  since this is $1$-displacement convex as discussed in \cite{CT}, and is minimized by $v^{(\infty)}\dd x$. 
Thus taking $\dd \mu = v^{(\infty)}(x)\dd x$ and $\dd \nu = v(x)\dd x$ in ${\mathcal M}_M$,
 (\ref{abovetan2}) specializes to
\begin{equation}\label{tala}
 {\rm W}_2^2(\mu,\nu)  \leq 
{\mathcal H}[v|v^{(\infty)}]\ .
\end{equation}

\noindent{\bf Proof of Lemma~\ref{entsecmo}:} Because of (\ref{tala}), it remains to show that 
$$\alpha[v|v^{(\infty)}] \leq 2\sqrt{\alpha[v^{(\infty)}]}{\rm W}_2(v,v^{(\infty)}) + {\rm W}_2^2(v,v^{(\infty)})\ .$$
To see this, let $\gamma$ be an optimal coupling of $v(x)\dd x$ and $v^{(\infty)}(x)\dd x$, and note that
$\sqrt{\alpha[v]} = \|x\|_{L^2(\R^2, \dd \gamma)}$ and  $\sqrt{\alpha[v^{(\infty)}]} = \|y\|_{L^2(\R^2, \dd \gamma)}$
where we write $x$ and $y$ to denote the functions $(x,y)\mapsto x$ and  $(x,y)\mapsto y$ respectively. Then 
by the triangle inequality,
\begin{equation}\label{alpdif}
|\sqrt{\alpha[v]} - \sqrt{\alpha[v^{(\infty)}]}| \leq  \|x-y\|_{L^2(\R^2, \dd \gamma)} =  {\rm W}_2(v,v^{(\infty)})\ .
\end{equation}
\qed

\section{Construction of properly dissipative weak solutions}

Define ${\mathcal K}_M$ to be the subset of ${\mathcal M}_N$ consisting of absolutely continuous measures $v(x)\dd x$ such that $E[v] < \infty$. 

We now implement the JKO scheme \cite{AGS,JKO}: Given a small time step $\tau >0,$  and given an initial datum $v_0 \in{\mathcal K}_M$, we recursively define a sequence $(v_n^{\tau})_{n \in \mathbb{N}}$    by $v_0^{\tau} = v_0$ and
\begin{equation}\label{discvp}
v_{n+1}^{\tau} = \mathrm{argmin}_{v \in {\mathcal K}_M} \left \{  \tau E[v|v^{(\infty)}] + \frac{1}{2} W_2^2(v,v_{n}^{\tau}) \right \},
\end{equation}

The first step of our analysis is to show that the variational scheme in (\ref{discvp}) has a solution and to derive the Euler-Lagrange equation for this variational scheme. 

\begin{lm}\label{ot}
Given $\tau>0,$ let $v_0 \in   {\mathcal  K}$ . 

\begin{itemize}
\item[(i)] The discrete  variational  scheme (\ref{discvp}) admits a solution $\{v_n^{\tau}\}_{n \in \mathbb{N}}$.

\item[(ii)] Each $v_n^{\tau}$ satisfies
\begin{equation}\label{reg3der}
(v_n^{\tau})_{xxx}  \in L_{loc}^{\infty}(\R),
\end{equation}
and
\begin{equation}\label{dervanish}
v_n^{\tau}(a) = 0   \implies (v_n^{\tau})_x(a) = 0, \quad a \in \R.
\end{equation}

\item[(iii)]  Let $\psi_n'$ be the the optimal transportation plan such that  $v_{n-1}^{\tau}(x)\dd x = (\psi_n')\#
v_{n}^{\tau}(x)\dd x$. Then
\begin{equation}\label{EL1}
\psi_n'(x) = x + \tau\left[ xv_n^\tau(x)  - v_n^\tau(x) (v_n^\tau)_{xxx}(x)\right]\ .
\end{equation}
\end{itemize}
\end{lm}

\medskip

The proof of this lemma is closely patterened on a proof of Otto \cite{O}. We present some details for the convenience of the reader since
Otto's energy functional differs from our in having the term $(1/2)\int_{\R}x^2v(x)\dd x$ replaced by the Lebesgue measure of the set 
$\{x\:\ v(x) > 0\}$.  The reader familiar with \cite{O}, or even \cite{JKO}, may wish to skip ahead to Lemma~\ref{Hdec}. For other readers, we point out
that (\ref{EL1}) is the Euler-Lagrange equation for the discrete time variational principle, and through (\ref{EL1}), one makes a direct 
connection with the thin film equation (\ref{veq}).

\noindent{\bf Proof of Lemma~\ref{ot}:}  
\emph{Proof of (i):}
 It is enough to show that for given $v_0 \in {\mathcal K}_M,$ there exists a solution of
\begin{equation}\label{ilki}
 v_1^{\tau} = \mathrm{argmin}_{v \in {\mathcal K}_M}  \left\{\tau E[v|v^{(\infty)}] + \frac{1}{2} W_2^2(v_0,v)\right\}.
\end{equation}

Let ${v^k}_{k \in \N} $ be a minimizing sequence in (\ref{ilki}). As $\{E(v^k|v^{(\infty)})\}_{k\in \N}$ is bounded, there exists $v_1 : \mathbb{R}  \to [0,\infty)$ such that

\begin{equation}\label{convk}
v^k  \to v_1, \quad \text{locally uniformly on} \quad \R
\end{equation}
for a subsequence and that
\begin{equation}\label{enbdvk}
E(v_1|v^{(\infty)})  \leq \liminf_{k \to \infty} E(v^k|v^{(\infty)}).
\end{equation}
Indeed, by (\ref{infbnd}) and the elementary bound
${\displaystyle 
 \sup_{x_1,x_2 \in \R, x_1 \neq x_2} \frac{v(x_1)-v(x_2)}{\sqrt{|x_1-x_2|}} \leq \sqrt{2E[v]}}$,
 the Arzela-Ascoli theorem may be used to prove (\ref{convk}). Due to (\ref{convk}) and  the boundedness of $\{(v^k)'\}_{k \in \N}$ in $L^2(\R)$ we have $v_1' \in L^2(\R)$,
$ (v^k)_x \rightharpoonup (v_1)_x,$ in $L^2(\R)$,
and in particular
\begin{equation}\label{l2derbd}
\int_{\R} |(v_1)_x|^2 \, \dd x \leq \liminf_{k \to \infty} \int_{\R} |(v^k)_x|^2 \, \dd x.
\end{equation}
Using the inequality (\ref{alpdif}),
for all $v_0,v_1 \in {\mathcal K}$,   we deduce that the boundedness of $\{W_2(v_0,v^k)\}_{k \in \N}$ implies that
$\left\{\int_{\R} x^2 v^k \, \dd x\right\}_{k \in \N}$  is bounded.
Then, by this,  (\ref{convk})  and  Fatou's lemma, we deduce
$$ \int_{\R} x^2 v_1(x) \, \dd x < \infty.$$

Also, $\int_{-r}^r v^k(x) \, \dd x   = M - \int_{\R \setminus [-r,r]} v^k \, \dd x \geq M - \frac{1}{r^2} \int_{\R} x^2  v^k \, \dd x$ and so $\lim_{r \uparrow \infty} \int_{-r}^r v^k \, \dd x =M$, uniformly in $k.$ And thus, $\int_{\R} v_1 \, \dd x =1.$
Combining, we have $ E(v_1|v^{(\infty)}) \leq \liminf_{k \to \infty} E(v^k|v^{(\infty)})$
along a subsequence. Next, the fact that
\begin{equation}\label{w2bdd}
W_2^2(v_0,v_1) \leq \liminf_{k \to \infty} W_2^2(v_0,v^k).
\end{equation}
follows from the lower semicontinuity of the Wasserstein distance, see \cite{O} or \cite{V}. This concludes the proof of \emph{(i)}.

\emph{Proof of (ii) and (iii)} The key step here is a variational argument of Otto \cite{O} showing that, in our case, at each $n$
the optimal transportation plan $\psi_n'$ such that  $v_{n-1}^{\tau}(x)\dd x = (\psi_n')\#
v_{n}^{\tau}(x)\dd x$ is related to $v_{n}^{\tau}$ through
\begin{equation}\label{1stvar}
\frac{1}{\tau} \int_{\R} (x-\psi_n'(x)) \xi(x) v_n^{\tau}(x) \, \dd x + \int_{\R} \left[ x \xi(x) v_n^{\tau}(x) - \frac{3}{2} (v_n^{\tau})_x^2 \xi'(x)- v_n^{\tau} (v_n^{\tau})_x \xi''(x)  \right] \, \dd x = 0.
\end{equation}

To make the variation of $v_n^{\tau}$, consider any  $\xi \in C_0^{\infty}(\R)$,  and
 define $\tilde{v}(x)\dd x = (Id+\eps \xi)\#(v_n^{\tau}(x)\dd x)$. 
 We now work out the effects of this variation on each term in the functional being minimized in (\ref{ilki}), starting with the Wasserstein distance. 
 
 Let $\gamma$ denote the optimal coupling of 
 $v_{n}^{\tau}(x)\dd x$ and  $v_{n-1}^{\tau}(x)\dd x$. 
 Then $\tilde{\gamma} := ((Id + \eps \xi) \otimes Id)\# \gamma$ is some (non-optimal) 
 coupling of  $\tilde{v}\dd x$ and  $v_{n-1}^{\tau}(x)\dd x$, and hence

\begin{equation}\label{ineqwas}
\frac{1}{2} W_2^2(\tilde{v},v_{n-1}^{\tau}) - \frac{1}{2} W_2^2(v_{n-1}^{\tau},v_{n}^{\tau}) \leq \eps \int_{\R} (x-\psi'(x)) \xi(x) v_{n}^{\tau} \, \dd x + O(\eps^2).
\end{equation}

As for the second moments, we have
\begin{equation}\label{es2mom}
\frac{1}{2} \int_{\R} x^2 \tilde{v}(x) \, \dd x - \frac{1}{2} \int_{\R} x^2 v_{n}^{\tau} (x) \, \dd x = \eps \int_{\R} x \xi(x) v_{n}^{\tau} (x) \, \dd x + O(\eps^2).
\end{equation}

Next, since
\begin{equation}\label{cofvar}
\tilde{v}(y) = v_{n}^{\tau}(x(y)) \left|\frac{\dd x}{\dd y}\right| \quad \text{with} \quad \frac{\dd x}{\dd y} = 1-\eps \xi'(y) + O(\eps^2)\ ,
\end{equation}
\begin{equation}\label{esenrgy}
\frac{1}{2}\int_{\R} \tilde{v}_x^2 \, \dd x - \frac{1}{2} \int_{\R}(v_{n}^{\tau})_x^2 \, \dd x = - \eps \int_{\R} \left[ \frac{3}{2} (v_{n}^{\tau})_x^2 \xi'(x) + v_{n}^{\tau} (v_{n}^{\tau})_x \xi''(x)   \right] \, \dd x + O(\eps^2)
\end{equation}
Now  (\ref{1stvar}) follows by combining (\ref{ineqwas}),(\ref{es2mom}),(\ref{esenrgy}) and replacing $\xi$ by $-\xi$.

Next, we complete the proof of {\em (ii)}: 
As in  Proposition 1.3 of \cite{O}, because any $v\in {\mathcal K}$ is continuous, the set  $\{ v_n^{\tau}>0\}$ consists of countably many intervals. Let $I = (a,b)$ be such an interval.
Then we have from  (\ref{1stvar}) that for some $f \in L^{\infty}_{loc}(\R)$ (expressible n terms of $\psi_n'$), whenever $\xi$ that is smooth and compactly supported in $I$,   
$$\int_{\R}\left[ \frac{3}{2} (v_n^{\tau})_x^2 \xi'(x)+ v_n^{\tau} (v_n^{\tau})_x \xi''(x)\right] \dd x  =  \int_{\R}\xi(x) f(x)v_n^{\tau}(x)\dd x\ .$$
Hence, taking distributional derivatives, 
\begin{equation}\label{dist}
\left[  \frac{3}{2}(v_n^{\tau})_x^2 - (v_n^{\tau}(v_n^{\tau})_x)_x   \right]_x = v_n^{\tau} f, \quad \text{in} \quad D^*(I)\ .
\end{equation}
It follows that $v_n^{\tau} \in H_{loc}^{3,1}(I)$, see \cite{O} for details.  For any $v\in  H_{loc}^{3,1}(I)$,
\begin{equation}\label{dist2}\frac{3}{2}(v_x)^2_x = 3v_xv_{xx}\qquad{\rm and}\qquad
(vv_x)_{xx} = 3v_xv_{xx} + vv_{xxx}\ .
\end{equation}
Hence we can rewrite (\ref{dist}) as
\begin{equation}\label{dist3}
 (v_n^{\tau})_{xxx}(x)  =  f(x)\quad \text{on} \quad I\ ,
 \end{equation}
 for the same $f$. 
Next, as noted above, $\{ v_n^{\tau}>0\}$ is a union of countably many  intervals $\{(a_j, b_j)\}_{j \geq 1}$. Integrating by parts, using the fact that
$v_n^{\tau}$ vanishes on the endpoints of each interval, 
\begin{eqnarray}
-\int_{\R} \xi v_n^{\tau} f \, \dd x &=& \int_{\R} \left[ -\frac{3}{2}(v_n^{\tau})_x^2 \xi' - v_n^{\tau} (v_n^{\tau})_x \xi'' \right] \, \dd x \nonumber \\
 		    &=& \sum_{j=1}^{\infty} \int_{a_j}^{b_j}  \left[ -\frac{3}{2}(v_n^{\tau})_x^2 \xi' - v_n^{\tau} (v_n^{\tau})_x \xi'' \right] \, \dd x \nonumber \\
	    &=& \sum_{j=1}^{\infty} \left[  -\frac{1}{2}(v_n^{\tau})_x^2 \xi |_{a_j}^{b_j} - \int_{a_j}^{b_j} v_n^{\tau} (v_n^{\tau})_{xxx} \xi \, \dd x \right] \nonumber\\ 	         &=& \sum_{j=1}^{\infty} -\frac{1}{2} (v_n^{\tau})_x^2 \xi |_{a_j}^{b_j} - \int_{a_j}^{b_j} \xi v_n^{\tau} f, \quad \forall \xi \in C_0^{\infty}(\R).
\end{eqnarray}
Equivalently,
$$ - \frac{1}{2} \sum_{j=1}^{\infty} (v_n^{\tau})_x^2 \xi |_{a_j}^{b_j} = 0, \quad \forall \xi \in C_0^{\infty}(\R).$$
Since $\xi$ is an arbitrary smooth compactly supported function, this implies  that  $ (v_n^{\tau})_x(a_j) = 0 = (v_n^{\tau})_x(b_j),$ for all $j \geq 1.$

Finally, we have enough regularity to deduce {\em (iii)} from (\ref{1stvar}), which implies that
$$\psi_n'(x) = x + \tau\left[ xv_n^\tau(x) + \frac{3}{2}(v_n^\tau)^2_{xx}(x) - (v_n^\tau(v_n^\tau)_x)_{xx}(x)\right]\ .$$
Simplifying this using (\ref{dist2}),  we obtain (\ref{EL1}). 
\qed

At this stage, we depart from Otto's analysis in \cite{O}. Our next goal is to show that the entropy $H[v_n^{\tau}|v^{(\infty)}]$
is monotone decreasing in $n$ and to relate the decrease to what what one would guess by a formal differentiation argument with the continuous time evolution equation. For this we  exploit the displacement convexity of the functional $H[v|v^{(\infty)}]$. 
First, we make a definition:

\begin{defi}[Entropy dissipation]  The functional $D$ is defined on ${\mathcal K}_M$ by
\begin{equation}\label{disdef}
D[v] :=  \int_{\mathbb{R}} \left( x + \sqrt{6}v_x(x) \right)^2 v(x) \, \dd x\  .
\end{equation}
\end{defi}

\begin{lm}\label{Hdec}
Fix $n  \in \mathbb{N}$ and let $v_n^{\tau}$ be a solution of  (\ref{discvp}). Then, the following inequality holds
\begin{eqnarray}\label{Hdiscvx}
\frac{H[v_n^{\tau}|v^{(\infty)}] - H[v_{n-1}^{\tau}|v^{(\infty)}]  }{\tau} &\leq& - D[v_n^{\tau}] -
\frac{\sqrt{6}}{24}  \int_{\R} (v_n^{\tau})^{-3/2} (v_n^{\tau})_x^4\nonumber\\
&\leq& -2 H[v_n^{\tau}|v^{(\infty)}]
\end{eqnarray}

\end{lm}

\noindent{\bf Proof:}  The second inequality follows from the first using the entropy-entropy dissipation inequality 
$2H[v|v^{(\infty)}] \leq  D[v]$, as explained in \cite{CT}. Hence we must prove the first inequality. 

Let $(v_n^{\tau})^{(t)}$, $0\leq t \leq 1$, denote the displacement interpolation between  $v_n^{\tau}$ and $v_{n-1}^{\tau}$. Since $H[v]$
is a displacement convex functional, it follows that for $t\in (0,1)$
$$H[v_{n-1}^{\tau}] - H[v_{n}^{\tau}] \geq \frac{1}{t}\left(H[(v_n^{\tau})^{(t)}] - H[v_{n}^{\tau}]\right)\ ,$$
and moreover, the right hand side is monotone decreasing as $t$ tends to zero .

By a standard computation, \cite{AGS,V},
$$\lim_{t \to 0} \frac{1}{t}\left(H[(v_n^{\tau})^{(t)}|v^{\infty}] - H[v_n^{\tau}|v^{(\infty)}]\right) \geq -\int_\R \left( \frac{\delta H}{\delta v}[v_n^{\tau}]\right)
(\psi_n'(x) -x)_x \, \dd x\ ,$$
where
$$ \frac{\delta H}{\delta v}[v_n^{\tau}](x)  = \frac{x^2}{2} + \sqrt{6}  \sqrt{v_n^{\tau}}(x)\ .$$

Integrating by parts and using the Euler-Lagrange equation (\ref{EL1}), 
\begin{eqnarray}\label{disconineq}
H[v_{n-1}^{\tau}|v^{\infty}] - H[v_n^{\tau}|v^{(\infty)}] &\geq& \int_{\R} [\frac{x^2}{2}+\sqrt{6}\sqrt{v_n^{\tau}}]_x (\psi_n'(x)-x) v_n^{\tau} \, \dd x \nonumber \\
				  & =& \tau  \int_{\R} [x+\sqrt{6}(\sqrt{v_n^{\tau}})_x] ( x -(v_n^{\tau})_{xxx} ) v_n^{\tau} \, \dd x, 	
\end{eqnarray}
where we have used (\ref{EL1}). Adding and subtracting $x+\sqrt{6}(\sqrt{v_n^{\tau}})_x,$ we deduce that
\begin{equation}\label{finineqdc}
H[v_{n-1}^{\tau}|v^{(\infty)}]- H[v_n^{\tau}|v^{(\infty)}] \geq  \tau \int_{\R} \left[ x + \sqrt{6}(\sqrt{v_n^{\tau}})_x \right]^2 v_n^{\tau} \, \dd x  + J
\end{equation}
where
\begin{equation}\label{finineqdc2}
J := - \tau \int_{\R} \left[ x + \sqrt{6} (\sqrt{v_n^{\tau}})_x \right] \left[ \sqrt{6} (\sqrt{v_n^{\tau}})_x + (v_n^{\tau})_{xxx} \right] v_n^{\tau} \, \dd x.
\end{equation}
Integrating by parts, using Lemma~\ref{regresults}, proved just below,  to justify certain of these integration by parts, we obtain
\begin{eqnarray}\label{J+}
 J = \frac{\sqrt{6}}{3} \tau \int_{\R} (v_n^{\tau})^{3/2} \, \dd x &+& \frac{\sqrt{6}}{2} \tau \int_{\R} (v_n^{\tau})^{1/2}(v_n^{\tau})_{xx}^2 \, \dd x \nonumber \\
&+& \frac{\sqrt{6}}{24} \tau \int_{\R} (v_n^{\tau})^{-3/2} (v_n^{\tau})_x^4 \, \dd x\ ,
\end{eqnarray}
each term of which is non-negative.
Combining (\ref{J+}) and (\ref{finineqdc}) we obtain the result.
\qed

\begin{lm}\label{regresults}
Let $v_n^{\tau}$ be the $n$th step in a  solution of (\ref{discvp}). Then,
\begin{equation}\label{goodibp}
   \frac{(v_n^{\tau})_x^3}{(v_n^{\tau})^{1/2}}  \in L_{loc}^{\infty}(\R), \quad \text{and} \quad \frac{(v_n^{\tau})_x^3(a)}{(v_n^{\tau})^{1/2}(a)} \to 0 \quad \text{if} \quad v_n^{\tau}(a) = 0.
\end{equation}
\end{lm}

\noindent{\bf Proof:} 
Without loss of generality assume that $v_n^{\tau}(0) = 0 = (v_n^{\tau})_x(0).$ Assume on the contrary  that  $\frac{(v_n^{\tau})_x^3}{(v_n^{\tau})^{1/2}} \geq \delta^3 > 0$ for some $\delta >0.$ This implies $(v_n^{\tau}(x))^{1/2} \geq C_1 x^{3/5}$ for some constant $C_1>0.$ Indeed, we easily deduce from the assumption that $\frac{(v_n^{\tau})_x}{(v_n^{\tau})}\geq \delta > 0$ which implies easily that  $[(v_n^{\tau})^{5/6}]_x \geq C \delta.$ Then, $(v_n^{\tau}(x))^{5/6} \geq C \delta x$ with $C>0$ holds. This immediately implies the assertion.  Now, by (\ref{reg3der}) we get $(v_n^{\tau})_x \leq C_2 x,$ for some constant $C_2>0.$ Hence, $ \frac{(v_n^{\tau})_x^3}{(v_n^{\tau})^{1/2}} \leq C \frac{x^3}{x^{3/5}}  = C x^{12/5} \to 0$  as $x \to 0.$ This is a contradiction.
\qed

We next examine the behavior of second moments along the discrete scheme.

\begin{lm}\label{displalph}
Fix $N, M \in \mathbb{N}$ with $N \geq M$ and let $v_n^{\tau}$ be a solution of (\ref{discvp}). Then, the following inequality holds true.
\begin{eqnarray}\label{alphaineq1}
\alpha[v_N^{\tau}|v^{(\infty)}] - \alpha[v_M^{\tau}| v^{(\infty)}]     &+& 5 \tau \sum_{j = M+1}^N  \alpha[v_j^{\tau}|v^{(\infty)}] \nonumber \\
&\geq & 3\tau  \sum_{j= M+1}^N E[v_j^{\tau}|v^{(\infty)}] - \tau  E[v_0|v^{(\infty)}] \ .
\end{eqnarray}
\end{lm}

\noindent{\bf Proof:} Since
$$\alpha[v_n^{\tau}|v^{(\infty)}] -\alpha[v_{n-1}^{\tau}|v^{(\infty)}] = \alpha[v_n^{\tau}] -\alpha[v_{n-1}^{\tau}] = \int_{\R}(|x|^2 - |\psi_n'(x)|^2)
v_n^{\tau}(x)\dd x\ ,$$
a simple computation using the
Euler-Lagrange equation (\ref{EL1}) and the regularity results to integrate by parts, one obtains
\begin{equation}\label{barcba}
\alpha[v_n^{\tau}|v^{(\infty)}] -\alpha[v_{n-1}^{\tau}|v^{(\infty)}] \geq -2 \tau\alpha[v_{n-1}^{\tau}|v^{(\infty)}] + 3\tau\beta[v_{n-1}^{\tau}|v^{\infty}] - \frac{1}{2}{\rm W}_2^2(v_n^{\tau},v_{n-1}^{\tau}).
\end{equation}
Thus, 
\begin{equation}\label{estvgood}
\alpha[v_N^{\tau}|v^{(\infty)}] - \alpha[v_{M}^{\tau}|v^{(\infty)}] + 5 \tau \sum_{j=M+1}^N \alpha[v_j^{\tau}|v^{(\infty)}] \geq 3 \tau \sum_{j=M+1}^N E[v_j^{\tau}|v^{(\infty)}] -\frac{1}{2} \sum_{j=M+1}^N {\rm W}_2^2(v_{j-1}^{\tau},v_j^{\tau}),
\end{equation}
which is obtained by summing up the estimate (\ref{barcba}).

Now, observe that thanks to the variational structure of (\ref{discvp}), we obtain the following estimate for free
\begin{equation}
\sum_{k=1}^N {\rm W}_2^2(v_{k-1}^{\tau},v_k^{\tau}) \leq 2 \tau \Big( E[v_0|v^{(\infty)}] - E[v_N^{\tau}|v^{(\infty)}]    \Big),
\end{equation}
and hence using this in (\ref{estvgood}) we obtain (\ref{alphaineq1}).
\qed

We next control the  fourth moments.

\begin{lm}\label{disc4mbd}
Fix $n \in \mathbb{N}$ and let $v_n^{\tau}$ be a solution of (\ref{discvp}) with initial data $v_0\in {\mathcal K}$
having a finite fourth moment; i.e.,  $M_4(v_0) < \infty$. Then $M_4(v_n^{\tau})$ is bounded uniformly in $n$. Indeed,
there is an explicit constant $K_3$ depending only on $M$ and $E[v_0]$ such that 
\begin{equation}\label{disc4mombdd}
M_4(v_n^{\tau})   \leq  M_4(v_0) + K_3H[v_0|v^{(\infty)}]\ .
\end{equation}
\end{lm}

\noindent{\bf Proof:} 
By the displacement convexity of the functional $M_4(f)$ we deduce that
\begin{eqnarray}\label{4mtdiscbd}
M_4(v_{n-1}^{\tau}) &\geq&  
 M_4(v_n^{\tau}) + \frac{d}{d\tau}\int_{\R} (x+ \tau (x-(v_n^{\tau})_{xxx}))^4 v_{n}^{\tau} \, \dd x \bigg|_{\tau = 0} \nonumber \\
&=&  M_4(v_n^{\tau}) + 4 \tau  \int_{\R} x^3  (x-(v_n^{\tau})_{xxx})) v_{n}^{\tau} \, \dd x\ \nonumber\\
&=&  M_4(v_n^{\tau}) + 4 \tau  M_4(v_n^{\tau})  -\tau \int_{\R} x^3(v_n^{\tau})_{xxx} v_{n}^{\tau} \, \dd x\  .\nonumber\\
\end{eqnarray}
We now integrate by part on the last term, using Lemma \ref{regresults} to  justify the calculations. We obtain
\begin{equation}\label{4momfin1}
\frac{M_4(v_{n}^{\tau}) - M_4(v_{n-1}^{\tau})  }{\tau} \leq - 4 M_4(v_n^{\tau}) - 12\int_{\R} (v_n^{\tau})^2 \, \dd x + 18 \int_{\R} x^2 [(v_n^{\tau})_x]^2 \, \dd x.
\end{equation}
The last term on the right hand side of (\ref{4momfin1}) is estimated as follows:
\begin{eqnarray}
\int_{\R} x^2 [(v_n^{\tau})_x]^2 \, \dd x &=& \int_{\R} x^2 (v_n^{\tau})^{3/4}(v_n^{\tau})^{-3/4} [(v_n^{\tau})_x]^2 \, \dd x\nonumber\\
&\leq&\|v_n^{\tau}\|_\infty^{3/4}
  (M_4(v_n^{\tau}))^{1/2}\left(\int_{\R} (v_n^{\tau})^{-3/2}(v_n^{\tau})_x^4 \, \dd x\right)^{1/2}\ ,  \nonumber\\
  &\leq& 4M_4(v_n^{\tau}) + \frac{\|v_n^{\tau}\|_\infty^{3/2}}{16} \int_{\R} (v_n^{\tau})^{-3/2}(v_n^{\tau})_x^4 \, \dd x\ .\nonumber
  \end{eqnarray}
Combining this with (\ref{4momfin1}) and using (\ref{infbnd}) and Lemma~\ref{Hdec}, we obtain
$$
M_4(v_{n}^{\tau}) - M_4(v_{n-1}^{\tau})    \leq K_3
(H[v_{n-1}^{\tau}|v^{(\infty)}] - H[v_{n}^{\tau}|v^{(\infty)}]) \ .$$
Telescoping the sums gives the result. \qed

We are now ready to prove Theorem~\ref{weakex},   the existence of properly dissipative weak solutions:

\medskip
\noindent{\bf Proof of Theorem~\ref{weakex}:} Define $v^{(\tau)}(x,t) = v_n^{\tau}(x)$ for $n\tau \leq t \leq (n+1)\tau$. It is then standard to
 show \cite{AGS,O,JKO} that
from the family $\{v^{(\tau)}\}_{\tau>0}$, one can extract a weakly convergent subsequence, and that the weak limit is a weak solution of
(\ref{veq}) in the sense of (\ref{int-eqn}). See \cite{O} for such an argument.

Next, for any $\tau>0$ and $t= n\tau$, we have from Lemma\ref{Hdec} that
$$2\int_0^t H(v|v^{(\infty)}) \, dt + H[v(\cdot,t)|v^{(\infty)}] \leq H[v(\cdot,0)|v^{(\infty)}]\ .$$
A standard convexity and lower semicontinuity argument shows that this inequality is preserved along weakly converging subsequences. This 
proves (\ref{H-dissp}).  The argument is relatively straightforward since the ``small'' side of the inequality is a weakly lower semicontinuous function
of $v(x,t)$, while the ``large'' side of the inequality only depends on the initial data.

A more involved argument is required to prove (\ref{alph-disp}) since in this case the solution $v(x,t)$, and not only the initial data,
occurs on both sides of the inequality. We therefore need continuity, and not only lower semicontinuity, of the functionals on the
 ``large'' side. 

This is provided by the uniformly bounded fourth moment.  On account of this, we
 conclude that along any weakly convergent subsequence, $\{v^{(\tau_k)}\}$
$$\lim_{k\to\infty}\alpha[v^{(\tau_k)}(\cdot,t)] =\alpha[v(\cdot,t)]\ ,$$
 while again a standard convexity and lower semicontinuity argument
shows that 
$$\lim_{k\to\infty}E[v^{(\tau_k)}(\cdot,t)]  \geq E[v(\cdot,t)]\ ,$$
Since $E$ is decreasing along each solution of the discrete scheme, (\ref{alph-disp}) holds with $v^{(\tau_k)}$ in place of $v$, and then
by what we have said above, the inequality is preserved in the limit.  This proves  (\ref{alph-disp}).
\qed.

\noindent{\bf Proof of Theorem~\ref{weakpro};} Since $v$ is a properly dissipative weak solution is satisfies (\ref{H-dissp}), and hence
$$H[v(\cdot,t)|v^{(\infty)}]  \leq e^{-2t}H[v_0|v^{(\infty)}] \ .$$
Then by Lemma~\ref{entsecmo},
$$\alpha[v(\cdot,t)|v^{(\infty)}] \leq   e^{-t}2\sqrt{\alpha[v^{(\infty)}]}  \sqrt{  H[v_0|v^{(\infty)}]  } +  e^{-2t}H[v_0|v^{(\infty)}] \leq K e^{-t} \ ,$$
where $K$ depends  on $v_0$ only through $M$ and $H[v_0|v^{(\infty)}]$. Then since $v$ is 
a properly dissipative weak solution is satisfies (\ref{alph-disp}), and hence $v$ satsfies (\ref{expd3}).  This proves (\ref{conv1}).

Next, since $|x|^{2p} = (|x|^2)^{2-p}(|x|^4)^{p-1}$, H\"older's inequality yields
$$\int_{\R}|x|^{2p}|v-v^{(\infty)}|^2\dd x \leq \left(\int_{\R}|x|^{2}|v-v^{(\infty)}|^2\dd x\right)^{2-p}
\left(\int_{\R}|x|^{4}|v-v^{(\infty)}|^2\dd x\right)^{p-1}\ .$$
By the pointwise bound, there is a constant $C$ depending only on $E[v_0]$ such that $\|v-v^{(\infty)}\|_\infty \leq C$, and hence for another constant $C$ depending only on  $E[v_0]$ and $M_4[v_0]$, $\int_{\R}|x|^{4}|v-v^{(\infty)}|^2\dd x \leq C$. Now (\ref{conv2}) follows easily.
\qed

\section{Acknowledgments} \label{ack}

The work of E. Carlen is partially suppoted by N.S.F. grant DMS  0901632.
The work of S. Ulusoy is partially supported by N.S.F. grants DMS 0707949, DMS1008397 and FRG0757227. 
Both authors thank Univ. Paul Sabatier-Toulouse  and IPAM, UCLA for  hospitality during collaborative visits  when  part of this work was completed.
E. Carlen thanks the CMA  at the University of Oslo and S. Ulusoy thanks Rutgers University for hospitality during other 
collaborative visits.

\end{document}